\documentclass[12pt]{article}
\usepackage{latexsym}
\date{}
\newtheorem{lemma}{Lemma}[section]

\newtheorem{remark}{Remark}[section]

\makeatletter 
\@addtoreset{equation}{section}
\makeatother 
 
\title{Generalization of the  Borel-Cantelli Lemma}

\author{A. Stepanov\thanks{ alexei@step.koenig.ru; Department of Mathematics,
Kaliningrad State Technical University, Sovietsky Prospect 1,
Kaliningrad, 236000 Russia.}}

\date{\begin{abstract} In the present note a generalization of
Borel-Cantelli Lemma is derived.
\end{abstract}}

\begin{document}

\maketitle  \vspace{3mm}\noindent


\vspace{3mm}

\section{Introduction} Suppose in the following $A_1,A_2,\ldots$ is
a sequence of events on a common probability space and
$\overline{A}_i$ denotes the complement of the event $A_i$. The
Borel-Cantelli Lemma (Lemma \ref{lemma1}) is very important for
producing strong limit theorems.
\begin{lemma}\label{lemma1}
If for any sequence  $A_1,A_2,\ldots$ of events
\begin{equation}\label{1}
\sum_{n=1}^\infty P\{A_n\}<\infty,
\end{equation}
then $P\{A_n\ i.o.\}=0$. If $A_1,A_2,\ldots$ is a sequence of
independent events and if $\sum_{n=1}^\infty P\{A_n\}=\infty$,
then $P\{A_n\ i.o.\}=1$.
\end{lemma}
The independence condition in the second part of this lemma have
been weakened in some investigations. We mention here the works of
Chung and Erdos (1952), Erdos and Renyi (1959), Lamperti (1963),
Kochen and Stone (1964), Spitzer (1964), and Petrov (2002).

The first part of the Borel-Cantelli Lemma was generalized in
Barndorff-Nielsen (1961) (Lemma \ref{lemma2}).
\begin{lemma}\label{lemma2}
Let $A_n\ (n\geq 1)$ be a sequence of events satisfying
$P\{A_n\}\rightarrow 0$. Let also
\begin{equation}\label{2}
\sum_{n=1}^\infty P\{ \overline{A}_nA_{n+1}\}<\infty.
\end{equation}
Then $P\{A_n\ i.o.\}=0$.
\end{lemma}
The results of Lemma \ref{lemma2}  holds true if in (\ref{2}) the
events $\overline{A}_nA_{n+1}$ are substituted with the events
$A_n\overline{A}_{n+1}$. Observe that condition (\ref{2}) in the
Barndorff-Nelson Lemma is weaker than  condition (\ref{1}) in the
Borel-Cantelli Lemma.

In the present note we  propose further generalization of Lemma
\ref{lemma1}.
\section{Results}
\begin{lemma}\label{lemma3}
Let $A_n\ (n\geq 1)$ be a sequence of events satisfying
$P\{A_n\}\rightarrow 0$. Let also  for some $m\geq 0$
\begin{equation}\label{3}
\sum_{n=1}^\infty P\{
\overline{A}_n\overline{A}_{n+1}\ldots\overline{A}_{n+m-1}
A_{n+m}\}<\infty.
\end{equation}
Then $P\{A_n\ i.o.\}=0$.
\end{lemma}
Observe that  condition (\ref{3}) when $m\geq 2$ is weaker than
condition (\ref{2}) in Lemma \ref{lemma2}.
\begin{remark} It can be formally
shown that for   any sequence $A_1,A_2,\ldots$ of events the
equality $P\{A_n\ i.o.\}=\alpha\in[0,1]$ holds iff $\
\lim_{n\rightarrow\infty}\sum_{k=0}^\infty
P\{\overline{A}_n\ldots\overline{A}_{n+k-1}A_{n+k}\}=\alpha$.
\end{remark}
\section*{References}
\begin{description}
{\small \item Barndorff-Nielsen, O. (1961).\ On the rate of growth
of the partial maxima of a sequence of independent identically
distributed random variables, {\it Math. Scand.}, {\bf 9},
383--394. \item Chung, K.L. and Erdos, P. (1952).\ On the
application of the Borel-Cantelli lemma, {\it Trans.  Amer. Math.
Soc.}, {\bf 72}, 179--186. \item Erdos, P. and Renyi, A. (1959).\
On Cantor's series with convergent $\sum 1/q_n$, {\it Ann. Univ.
Sci. Budapest. Sect. Math.}, {\bf 2}, 93--109. \item Kochen, S.B.
and Stone, C.J. (1964).\ A note on the Borel-Cantelli lemma, {\it
Illinois J.Math.}, {\bf 8}, 248--251. \item Lamperti, J. (1963).\
Wiener's test and Markov chains, {\it J. Math. Anal. Appl.}, {\bf
6}, 58--66. \item Petrov, V.V. (2002).\ A note on the
Borel-Cantelli lemma, {\it Statistic $\&$ Probability Letters},
{\bf 58}, 283--286. \item Spitzer, F. (1964).\ {\it Principles of
Random Walk}. - Van Nostrand, Princeton.}
\end{description}
\end{document}